\documentclass[12pt]{amsart}
\usepackage{epsfig}

\theoremstyle{plain}
\newtheorem{theorem}{Theorem}[section]
\newtheorem{conjecture}{Conjecture}[section]
\theoremstyle{remark}
\newtheorem{remark}[theorem]{Remark}
\newcommand\lra{\longrightarrow}

\newcommand{\nz}{{\mathbb N}}
\newcommand{\pz}{{\mathbb P}}
\newcommand{\qz}{{\mathbb Q}}

\begin{document}

\title[Modular threefolds of level eight]
{Modular Calabi--Yau threefolds of level eight}

\author{S\l awomir Cynk}

\thanks{Partially supported by DFG Schwerpunktprogramm 1094
(Globale Methoden in der komplexen Geometrie) and KBN grant no. 1 P03A 008 28.}


\address{Instytut Matematyki\\Uniwersytetu Jagiello\'nskiego\\
  ul. Reymonta 4\\30--059 Krak\'ow\\Poland}
\curraddr{Institut f\"ur Mathematik, Universit\"at Hannover,
Welfengarten~1, D--30060 Hannover, Germany}
\email{s.cynk@im.uj.edu.pl}

\author{Christian Meyer}
\address{Institut f\"ur Mathematik\\ Johannes
  Gutenberg-Uni\-ver\-sit\"at\\
  Stau\-din\-ger\-weg 9\\D--55099 Mainz\\Germany}
\email{cm@mathematik.uni-mainz.de}

\maketitle

\section{The $L$--series of rigid Calabi--Yau threefolds}
\label{sec:ls}

If $\tilde X$ is a Calabi--Yau threefold defined over $\mathbb Q$, and $p$
is a good prime (i.e., a prime such that the reduction of $\tilde X$
mod $p$ is nonsingular) then the map
\[
\operatorname{Frob}_p^* : H_{\text{\'et}}^i (\tilde{X}, \mathbb Q_l) \longrightarrow
H_{\text{\'et}}^i(\tilde{X}, \mathbb Q_l)
\]
on $l$--adic cohomology induced by the geometric Frobenius morphism gives rise to
$l$--adic Galois representations
\[
\rho_{l,i} : \operatorname{Gal}(\overline{\mathbb Q}/\mathbb Q)
\longrightarrow \operatorname{GL}_{b^i}(\mathbb Q_l).
\]

If a Calabi--Yau threefold $\tilde X$ is \emph{rigid} (i.e., $h^{1,2}(\tilde X)=0$
or equivalently $b^{2}(\tilde X)=2$) then $\tilde X$ is expected to be \emph{modular},
i.e., the $L$--series of the (semi-simplification of the) two-dimensional Galois representation
$\rho_{l,3}$ associated with the middle cohomology $H_{\text{\'et}}^3(\tilde{X}, \mathbb Q_l)$
equals the $L$--series of a cusp form $f$ of weight 4 for $\Gamma_{0}(N)$. The
precise conjecture has been formulated by Saito and Yui in \cite{SY}.
For details and examples the reader is referred to \cite{Yui} or \cite{Meyer}.

There are many examples of pairs of Calabi--Yau threefolds with an isomorphism
between some pieces of their middle \'etale cohomologies and the appropriate
Galois representations. In particular, if we can attach modular forms to these
pieces then these modular forms will be the same.
If on the other hand we detect the same modular forms in the middle \'etale
cohomologies of two Calabi--Yau threefolds then this should have a geometrical
reason:

\begin{conjecture}
(The Tate conjecture, as formulated in \cite[Conj. 5.8]{Yui})
If two isomorphic two-dimensional Galois representations $\rho_1$, $\rho_2$ occur in the
\'etale cohomology of varieties $X_1$, $X_2$ defined over $\qz$, then there should
be a correspondence between the two varieties (i.e., an algebraic cycle on the product
of the two varieties) defined over $\qz$, which induces an isomorphism between $\rho_1$
and $\rho_2$.
\end{conjecture}

Following \cite{HSvGvS} we will call two Calabi--Yau threefolds defined over $\qz$
{\em relatives} if the same (weight four) modular form occurs in their $L$--series.
Finding a correspondence between two relatives is a highly non-trivial task. It can be
induced by a birational map defined over $\qz$ or more generally by a finite map between
the two threefolds but this does not have to be the case. If a correspondence is induced
by a birational map then by a result of Batyrev (\cite{Bat}) the two Calabi--Yau threefolds
must have the same Betti (and Hodge) numbers.

In this note we will deal with Calabi--Yau modular threefolds associated with the
unique normalized weight four newform for $\Gamma_0(8)$. It can be written as a product
\[
f(q) = \eta(2\tau)^4\eta(4\tau)^4
\]
where $\eta(\tau) = q^{\frac{1}{24}} \prod_{n\in\nz} (1-q^n)$, $q=e^{2\pi i\tau}$,
is the {\em Dedekind eta function}. Some of the first Fourier coefficients $a_p$ of
$f(q)=\sum_{n=1}^{\infty} a_n q^n$ are:
\[
\begin{array}{|r|r|r|r|r|r|r|r|r|r|r|r|}
\hline
\rule[-3mm]{0mm}{8mm} \phantom{00}p&\phantom{00}2&\phantom{00}3&\phantom{00}5&\phantom{00}7&
\phantom{0}11&\phantom{0}13&\phantom{0}17&\phantom{0}19&\phantom{0}23&\phantom{0}29&\phantom{0}73\\
\hline
\rule[-3mm]{0mm}{8mm} a_p&0&-4&-2&24&-44&22&50&44&-56&198&154\\
\hline
\end{array}
\]
Note that by twisting with Legendre symbols we obtain newforms of different levels,
for example 16, 64, 72, 144, 200, 392, 400. If these newforms occur in the $L$--series
of some Calabi--Yau threefolds then we expect correspondences between them which are
defined over some finite extension of $\qz$.

We will verify the Tate conjecture for rigid level 8 Calabi--Yau
threefolds. More precisely we will give explicit correspondences
between most known examples. Elliptic fibrations play a special role
in the construction of correspondences, among them two modular curves
from the Beauville list. This is not surprising as the Galois representation
associated to the weight four level 8 newform is constructed using
the self-fiber product of the modular curves $X_{\Gamma}$, for
$\Gamma=\Gamma_{1}(4)\cap\Gamma(2)$ or $\Gamma=\Gamma_{0}(8)\cap \Gamma_{1}(4)$
(cf. \cite{Deligne}).

\section{List of level 8 rigid Calabi--Yau threefolds}
\label{sec:list}

In this section we list all known examples of rigid Calabi--Yau threefolds
associated with the weight four level 8 newform $f$. We are using the notation
$T_{h^{1,1}}$. Examples with the same Hodge numbers are listed separately (with
different superscripts) if no {\em birational} correspondence between them is known.

\subsection{Type $T_{70}$, $h^{1,1} = 70$}
\label{subsec:Fermi}

Consider the affine threefold $Z$ given by the equation
\[
x+\frac{1}{x}+y+\frac{1}{y}+z+\frac{1}{z}+t+\frac{1}{t} = 0.
\]
It has been studied in detail in \cite{PS}. It fibres via projection
to $\pz^1(x)$ into so called {\em Fermi surfaces} and is therefore
called {\em Fermi threefold}.

There are various birational Calabi--Yau models of $Z$. The invariants are
\[
\chi = 140, \qquad h^{1,1} = 70, \qquad h^{2,1} = 0.
\]
Verrill (\cite{Verrill}) considers $Z$ as a double cover of the toric variety
associated with the root lattice $A_1^3$.

Changing the signs of $z$ and $t$
we can rewrite the equation for $Z$ as
\[
\frac{(x+y)(xy+1)}{xy} = \frac{(z+t)(zt+1)}{zt}.
\]
This way we obtain a birational equivalence between $Z$ and the self-fiber
product of the Beauville elliptic surface $Y_{\Gamma}$ with $\Gamma=\Gamma_0(8)\cap\Gamma_1(4)$.
In the literature a small resolution of the self-fiber product is denoted by $W_0(8)$.
The singular fibers in this case are of the following types:
\[
\begin{array}{cccc}
I_8 & I_2 & I_1 & I_1\\
I_8 & I_2 & I_1 & I_1
\end{array}
\]
Fiber products of elliptic fibrations were first studied in detail by Schoen
(cf. \cite{Schoen}). Examples and correspondences with other Calabi--Yau
threefolds are investigated in sections \ref{sec:ref} and \ref{sec:rig}.

By homogenizing the equation for $Z$ we find a birational model as a quintic in $\pz^4$
defined by the equation
\[
w^2(xyz+xyt+xzt+yzt) = xyzt(x+y+z+t).
\]
This quintic is birationally equivalent with the double covering of $\pz^3$ branched
along the union of five planes and a Cayley cubic
which can be given by the equation
\[
u^2 = xyzt(x+y+z+t)(xyz+xyt+xzt+yzt).
\]

\subsection{Type $T_{70}^1$, $h^{1,1} = 70$}

Let $X$ be the double covering of $\pz^3$ branched along the octic surface
\[
xyzt(x-y)(y-z)(z-t)(t-x)=0.
\]
It occurs as arrangement no. 2 in \cite{CynkMeyer} and as arrangement no. 1 in \cite{Meyer}.
Resolving the singularities of $X$ we obtain a Calabi--Yau threefold with invariants
\[
\chi = 140, \qquad h^{1,1} = 70, \qquad h^{2,1} = 0.
\]

\subsection{Type $T_{50}$, $h^{1,1} = 50$}
\label{subsec:T_50}

Let $V_1$ be the double covering of $\pz^3$ branched along the octic surface
\[
xyzt(x+y)(y+z)(x-y-z-t)(x+y-z+t) = 0.
\]
It occurs as arrangement no. 29 in \cite{CynkMeyer} and as arrangement no. 32 in \cite{Meyer}.
Resolving the singularities of $V_1$ we obtain a Calabi--Yau threefold with invariants
\[
\chi = 100, \qquad h^{1,1} = 50, \qquad h^{2,1} = 0.
\]
Let $V_2$ be the double covering of $\pz^3$ branched along the octic surface
\[
xyzt(x+y)(x-y+z)(x-y-t)(x+y-z-t) = 0.
\]
It occurs as arrangement no. 44 in \cite{CynkMeyer} and as arrangement no. 69 in \cite{Meyer}.
Resolving the singularities of $V_2$ we obtain a Calabi--Yau threefold with the same invariants
as those of $V_1$. In fact a birational correspondence between $V_1$ and $V_2$ is exhibited
in section \ref{sec:constr}.

\subsection{Type $T_{46}$, $h^{1,1} = 46$}

Let $X$ be the double covering of $\pz^3$ branched along the octic surface
\[
xyzt(x+y)(x-y+z)(y-z-t)(x+z-t) = 0.
\]
It occurs as arrangement no. 62 in \cite{CynkMeyer} and as arrangement no. 93 in \cite{Meyer}.
Resolving the singularities of $X$ we obtain a Calabi--Yau threefold with invariants
\[
\chi = 92, \qquad h^{1,1} = 46, \qquad h^{2,1} = 0.
\]

\subsection{Type $T_{44}$, $h^{1,1} = 44$}

Let $X$ be the double covering of $\pz^3$ branched along the octic surface
\[
(x-t)(x+t)(y-t)(y+t)(z-t)(z+t)(x+y+z+t)(x+y+z-t)=0.
\]
It occurs as arrangement no. 87 in \cite{CynkMeyer} and as arrangement no. 238 in \cite{Meyer}.
Resolving the singularities of $X$ we obtain a Calabi--Yau threefold with invariants
\[
\chi = 88, \qquad h^{1,1} = 44, \qquad h^{2,1} = 0.
\]
The projective coordinate change
\[
(x:y:z:t)\mapsto(-\tfrac{y+z}2+t:-x-\tfrac{y+z}2:-\tfrac{y+z}2-t:\tfrac{y-z}2)
\]
transforms the branch locus into the octic surface given by
\[
(x-y)(x+y)(y-z)(y+z)(z-t)(z+t)(t-x)(t+x) = 0.
\]
The projective coordinate change
\[
(x:y:z:t)\mapsto(x-t:y-z:y+z:x+t)
\]
transforms this equation into
\[
xyzt(x+y+z-t)(x+y-z+t)(x-y+z+t)(-x+y+z+t) = 0.
\]

\subsection{Type $T_{40}$, $h^{1,1} = 40$}

Let $X$ be the double covering of $\pz^3$ branched along the octic surface
\[
xyzt(x+y+z+t)(x+y-z-t)(y-z+t)(x+z-t) = 0.
\]
It occurs as arrangement no. 241 in \cite{Meyer}. Resolving the singularities
of $X$ we obtain a Calabi--Yau threefold with invariants
\[
\chi = 80, \qquad h^{1,1} = 40, \qquad h^{2,1} = 0.
\]

\subsection{Type $T_{40}^1$, $h^{1,1} = 40$}

Consider a crepant resolution of a fiber product of two elliptic fibrations
with the following types of singular fibers:
\[
\begin{array}{cccc}
I_2 & I_2 & I_4 & I_4\\
I_2 & I_2 & I_4 & I_4
\end{array}
\]
It is the self-fiber product of the Beauville elliptic surface $Y_{\Gamma}$
with $\Gamma=\Gamma_1(4)\cap\Gamma(2)$. In the literature the resolution is
denoted by $W_1(4)$. It is a Calabi--Yau threefold with invariants
 \[
\chi = 80, \qquad h^{1,1} = 40, \qquad h^{2,1} = 0.
\]

\subsection{Type $T_{40}^2$, $h^{1,1} = 40$}

Consider a crepant resolution of a fiber product of two elliptic fibrations
with the following types of singular fibers:
\[
\begin{array}{cccc}
I_2 & I_2 & I_4 & I_4\\
I_2 & I_2 & D_4^* & I_2\\
\end{array}
\]
It is a Calabi--Yau threefold with invariants
 \[
\chi = 80, \qquad h^{1,1} = 40, \qquad h^{2,1} = 0.
\]

\subsection{Type $T_{40}^3$, $h^{1,1} = 40$}

Consider the complete intersection $X$ of four quadrics in $\pz^7$
given by the equations
\begin{align*}
  u_1^2 & = x^2-y^2,\\
  u_2^2 & = y^2-z^2,\\
  u_3^2 & = z^2-t^2,\\
  u_4^2 & = t^2-x^2.
\end{align*}
Resolving the singularities of $X$ we obtain a Calabi--Yau threefold $\tilde{X}$ with invariants
\[
\chi = 80, \qquad h^{1,1} = 40, \qquad h^{2,1} = 0.
\]
This is explained in detail in section \ref{sec:constr}.

\subsection{Type $T_{36}$, $h^{1,1} = 36$}

Consider a crepant resolution of a fiber product of two elliptic fibrations
with the following types of singular fibers:
\[
\begin{array}{ccc}
I_2 & D_6^* & I_2\\
D_6^* & I_2 & I_2
\end{array}
\]
It is a Calabi--Yau threefold with invariants
 \[
\chi = 72, \qquad h^{1,1} = 36, \qquad h^{2,1} = 0.
\]

\subsection{Type $T_{32}$, $h^{1,1} = 32$}
\label{subsec:Nygaard}

Consider the complete intersection $X$ of four quadrics in $\pz^7$
given by the equations
\begin{align*}
2 y_0^2 & = + x_0^2 - x_1^2 - x_2^2 - x_3^2,\\
2 y_1^2 & = - x_0^2 + x_1^2 - x_2^2 - x_3^2,\\
2 y_2^2 & = - x_0^2 - x_1^2 + x_2^2 - x_3^2,\\
2 y_3^2 & = - x_0^2 - x_1^2 - x_2^2 + x_3^2.
\end{align*}
It has 64 ordinary nodes as only singularities. There exist projective small resolutions
of all the nodes. The invariants of a small resolution $\tilde{X}$ of $X$ are
\[
\chi = 64, \qquad h^{1,1} = 32, \qquad h^{2,1} = 0.
\]
The threefold $\tilde{X}$ has been studied in detail by Nygaard and van Geemen in
\cite{NygaardvanGeemen}.

\subsection{Type $T_{32}^1$, $h^{1,1} = 32$}

Consider a crepant resolution of a fiber product of two elliptic fibrations
with the following types of singular fibers:
\[
\begin{array}{cccc}
I_4 & I_4 & I_2 & I_2\\
I_2 & I_2 & D_4^* & I_2
\end{array}
\]
It is a Calabi--Yau threefold with invariants
 \[
\chi = 64, \qquad h^{1,1} = 32, \qquad h^{2,1} = 0.
\]

\subsection{Type $T_{32}^2$, $h^{1,1} = 32$}

Consider a crepant resolution of a fiber product of two elliptic fibrations
with the following types of singular fibers:
\[
\begin{array}{cccc}
I_2 & I_2 & I_4 & I_4\\
I_4 & I_4 & I_2 & I_2
\end{array}
\]
It is a Calabi--Yau threefold with invariants
 \[
\chi = 64, \qquad h^{1,1} = 32, \qquad h^{2,1} = 0.
\]

\subsection{Type $T_{28}$, $h^{1,1} = 28$}

Let $X$ be the double covering of $\pz^3$ branched along the octic surface
\[
(x^2+y^2+z^2-t^2)(x^2+y^2-z^2+t^2)(x^2-y^2+z^2+t^2)(-x^2+y^2+z^2+t^2) = 0.
\]
Resolving the singularities of $X$ we obtain a Calabi--Yau threefold with invariants
\[
\chi = 56, \qquad h^{1,1} = 28, \qquad h^{2,1} = 0.
\]
The threefold $X$ is investigated in detail in \cite{Meyer}.

\subsection{Type $T_{16}$, $h^{1,1} = 16$}

Consider the complete intersection threefold $X\subset \pz^5$
defined by the equations
\begin{align*}
x_0^2 + x_1^2 + x_2^2 + x_3^2 & = 4 x_4 x_5,\\
x_4^4 + x_5^4 & = 2 x_0 x_1 x_2 x_3.
\end{align*}
The singular locus of $X$ consists of 12 singularities with local equation $x^2+y^2+z^4+t^4=0$
and 32 ordinary nodes. There exist projective small resolutions of all singularities of $X$.
The invariants of a small resolution $\tilde{X}$ of $X$ are
\[
\chi(\tilde{X}) = 32, \qquad h^{2,1}(\tilde{X}) = 0, \qquad h^{1,1}(\tilde{X}) = 16.
\]
The threefold $\tilde{X}$ is investigated in detail in \cite{Meyer}.

\section{Geometrical constructions}
\label{sec:constr}

In this section we investigate several geometrical constructions leading
to correspondences between Calabi--Yau threefolds.

Consider the two double octics $V_1$ and $V_2$ from \ref{subsec:T_50} and
observe that the branch loci are projectively equivalent to respectively
\begin{align*}
V_1:\qquad & zt(x+y)(x-y)(x+z)(t+y)(t+z)(t+y+z)=0,\\
V_2:\qquad & xt(x+z)(x-z)(x+y)(t+y)(t+z)(t+y+z)=0.
\end{align*}
The Cremona transformation
\[
(x,y,z,t)\mapsto(yz,xy,xz,xt),
\]
which is a birational involution of $\mathbb P^{3}$, transforms one of them
to the other. Thus the two double octics are birationally equivalent, as
announced in \ref{subsec:T_50}.

Consider a homogeneous polynomial $F(x,y,z,t)$ of degree 4. The $8:1$ map
\[
\pz^4(1,1,1,1,4) \lra \pz^4(1,1,1,1,4),
\]
\[
(x:y:z:t:w) \mapsto (x^2:y^2:z^2:t^2:xyztw)
\]
induces a correspondence between the two double octics given by
\[
w^2 = F(x^2,y^2,z^2,t^2)
\]
and by
\[
w^2 = x y z t F(x,y,z,t).
\]
If $F$ is a product of linear polynomials then the first octic is a union of
four quadric surfaces and the second octic is an arrangement of eight planes.
If we take $F(x,y,z,t) = (x-y)(y-z)(z-t)(t-x)$ then
$F(x^2,y^2,z^2,t^2) = (x-y)(x+y)(y-z)(y+z)(z-t)(z+t)(t-x)(t+x)$, and we obtain
a correspondence
\[
T_{44} \stackrel{8:1}\longrightarrow T_{70}^1.
\]
Note that a correspondence between these varieties was already constructed
implicitly by Nygaard and van Geemen in \cite{NygaardvanGeemen}.

If we take $F(x,y,z,t) = (x+y+z-t)(x+y-z+t)(x-y+z+t)(-x+y+z+t)$ then we obtain
a correspondence
\[
T_{28} \stackrel{8:1}\longrightarrow T_{44}.
\]

Now assume that the double octic $X$ is given by
\[
w^2 = f_1(x,y,z,t) f_2(x,y,z,t) f_3(x,y,z,t) f_4(x,y,z,t)
\]
with quadratic homogeneous polynomials $f_i(x,y,z,t)$. There is an obvious $8:1$
correspondence between $X$ and the intersection $Y$ of four quadrics in $\pz^7$ given by
\begin{align*}
u_1^2 = f_1(x,y,z,t),\\
u_2^2 = f_2(x,y,z,t),\\
u_3^2 = f_3(x,y,z,t),\\
u_4^2 = f_4(x,y,z,t).
\end{align*}

For example, this construction produces a correspondence
\[
T_{32}\stackrel{8:1}\longrightarrow T_{28}.
\]

As a second example consider the following intersection $X$ of four quadrics in $\pz^7$:
\begin{align*}
  u_1^2 & = x^2-y^2,\\
  u_2^2 & = y^2-z^2,\\
  u_3^2 & = z^2-t^2,\\
  u_4^2 & = t^2-x^2.
\end{align*}
This example is listed in section \ref{sec:list} as type $T_{40}^3$. We find the correspondence
\[
T_{40}^3\stackrel{8:1}\longrightarrow T_{44}.
\]

Consider also the composed map $X\lra \pz^4_{[1,1,1,1,4]}$ defined by
\[
\phi(x:y:z:t:u_1:u_2:u_3:u_4) = (x^2:y^2:z^2:t^2:xyztu_1u_2u_3u_4).
\]
It induces a correspondence
\[
T_{40}^3\stackrel{64:1}\longrightarrow T_{70}^1.
\]

Now we will describe a Calabi--Yau resolution $\tilde{X}$ of the singularities of $X$
and compute its Euler and Hodge numbers. The threefold $X$ is an iterated double covering
of $\pz^3$. The branch divisors are pairs of planes $(P_{1},P_{2})$, $(P_{3},P_{4})$,
$(P_{5},P_{6})$, $(P_{7},P_{8})$, intersecting along the
lines $l_{i}=P_{2i-1}\cap P_{2i}$.  Singularities of $X$ correspond
to the lines $l_{i}$. The lines $l_{i}$ intersect in four points
$Q_{1}=l_{1}\cap l_{2}$, $Q_{2}=l_{2}\cap l_{3}$, $Q_{3}=l_{3}\cap l_{4}$, $Q_{4}=l_{1}\cap l_{4}$.
We first blow up $\pz^3$ in the points $Q_{i}$ and then the strict
transforms of the lines $l_{i}$. After these blow--ups the singular
locus consists of eight ordinary nodes at the points $(1:\pm 1:\pm 1:\pm 1)\in\pz^3$.
The intersection of $X$ with the hyperplanes $u_1 = \sqrt{-1} u_2$, $u_3 = \sqrt{-1} u_4$,
$x = z$ is a surface which contains four of the nodes and is smooth in
these points (for the four remaining nodes take $x=-z$ instead). Thus there
exist projective small resolutions.

We first compute the Euler characteristic of the singular model
$X$. The idea is to stratify $\pz^3$ by the number of points in the
fibers of the iterated double cover. The generic fiber with 16
elements corresponds to a point outside the planes.

Counting points we see that the Euler characteristic of the sum of
planes is $8\cdot3-28\cdot 2+8+3\cdot 12=12$ (we count 8 planes, then subtract 28
lines and finally take into account 8 points on three planes and 12 on four planes). 

The set of fibers with 16 elements is the complement of the planes,
so its Euler characteristic is $4-12=-8$.

Fibers with eight elements correspond to the points on the eight planes outside
the double lines and also to the lines $l_{i}$ but without multiple points.
The Euler characteristic of the sum of lines on one plane is
$7\cdot 2-3\cdot1-2\cdot6=-1$, hence the Euler characteristic of the complement of
the lines in one plane is $3+1=4$.
The line $l_{1}$ contains two fourfold and two threefold points, the
Euler characteristic of the complement of multiple points in $l_{1}$
is $1$. In total the Euler characteristic of the set of fibers with eight elements
is $8\cdot4+4(2-4)=24$.

There are exactly eight points whose fiber consists of one element.
The rest contains fibers with four elements, its Euler characteristic is $4-(-8)-24-8=-20$.
Finally the Euler characteristic of $X$ is
\[
\chi(X)=16(-8)+8\cdot 24-4\cdot 20+8=-8.
\]

Now we will study the effect of blowing up a fourfold point.
A generic point of the exceptional divisor (outside the double lines) has
a fiber with 16 elements, the points on the lines have a fiber with 8 elements.
There are six double points: four of them have fibers with 4 elements and
two with 8 elements. Altogether the iterated cover of the exceptional
plane has Euler characteristic $16\cdot1-8\cdot2+4\cdot4=16$.
On $X$ there are four points over a fourfold point so the effect of blowing up fourfold
points is $4(16-4)=48$.  

Counting in the same manner for the blow--up of the double lines we obtain
that the set of fibers with 16 and 8 elements has Euler characteristic 0,
and there are four points with fiber with 4 elements. Thus the effect of blowing up
the double lines is $4(4\cdot4-4\cdot2)=32$. The small resolution of the eight
nodes increases the Euler characteristic by 8. Finally we compute
\[
\chi(\tilde{X}) = -8+48+32+8 = 80.
\]

Now it is possible to compute the Hodge numbers of $\tilde{X}$ with van Geemen's
point counting method. This requires counting points on the reduction mod $p$
for a sufficiently large prime $p$ such that Frobenius acts by multiplication
with $p$ on $H^2(\tilde{X})$. For details the reader is referred to \cite{Meyer}.
We find
\[
h^{1,1}(\tilde{X}) = 40, \qquad h^{2,1}(\tilde{X}) = 0,
\]
so $\tilde{X}$ is rigid.

\section{Rational elliptic fibrations}
\label{sec:ref}

Consider an arrangement $D$ of eight planes given by an equation
\[
f_{1}\cdot\ldots\cdot f_{8}=0,
\]
and assume that among the eight planes there are two disjoint quadruples
intersecting in a point each. After renumbering the equations and changing coordinates
we can assume that $f_{1},\dots,f_{4}$ depend only on $x,y,z$, whereas
$f_{5},\dots,f_{8}$ depend only on $y,z,t$. 

Let $S$ and $S'$ be the double coverings of $\pz^2$ branched along the
corresponding sums of four lines. Then in appropriate affine
coordinates (f.i., $z=1$) they can be written as follows:
\begin{align*}
S & = \{(x,y,u)\in \mathbb C^{3}:u^{2}=f_{1}(x,y,1)\cdot\ldots\cdot f_{4}(x,y,1)\}\\
S' & = \{(y,t,v)\in \mathbb C^{3}:v^{2}=f_{5}(y,1,t)\cdot\ldots\cdot f_{8}(y,1,t)\}
\end{align*}

This exhibits (birationally) both surfaces as elliptic fibrations. Moreover,
the map
\[
((x,y,u),(y,t,v)) \mapsto (x,y,t,uv)
\]
is a rational, generically $2:1$ map from their fiber product to
the double covering of $\pz^3$ branched along the octic surface $D$.

We will study elliptic fibrations coming from the above
construction. We can get the following sequences of singular fibers:
\[
(I_{2},I_{2},I_{4},I_{4}),\quad(I_{2},I_{2},D_{6}^{*}),
\quad (I_{2},I_{2},I_{2},D_{4}^{*}),\quad (I_{2},I_{2},I_{2},I_{2},I_{4}).
\]

The first two examples are unique whereas the other two change in a one
parameter family. In the table we give examples of explicit equations for
the branch locus of corresponding double quartic elliptic fibrations. We
also include the types and coordinates of the singular fibers and the Picard
number $\rho(E_{w})$ of the generic fiber (which can easily be computed using
the Zariski lemma).
\[
\def\arraystretch{1.3}
\def\arraycolsep{2mm}
\begin{array}[t]{|c|ccccc|c|}
 \hline
\def\multicolumn#1#2{}
 \quad S_{1} \quad&I_{2}&I_{2}&I_{4}&I_{4}&&{x(x+t)(x+z)(x+z+t)}\\
 \cline{2-6}
 &1&-1&0&\infty&&{\rho(E_{w})=1}\\
\hline 
 S_{2}&I_{2}&I_{2}&I_{4}&I_{4}&&{x(x+z+t)(x+z-t)(x+2z)}\\
 \cline{2-6}
 &0&\infty&-1&1&&{\rho(E_{w})=1}\\
\hline 
 S_{3}&I_{2}&I_{2}&I_{2}&D_{4}^{*}&&{x(x+t)(x+\lambda t)(x+z)}\\
 \cline{2-6}
 &0&1&\lambda&\infty&&{\rho(E_{w})=2}\\
\hline 
 S_{4}&I_{2}&I_{2}&I_{2}&D_{4}^{*}&&{t(x+\lambda z)(x+z)(x+\lambda t)}\\
 \cline{2-6}
 &0&1&\lambda&\infty&&{\rho(E_{w})=2}\\
\hline 
 S_{5}&I_{2}&I_{2}&D_{6}^{*}&&&{xt(x+z)(x+t)}\\
 \cline{2-6}
 &1&0&\infty&&&{\rho(E_{w})=1}\\
\hline 
 S_{6}&I_{2}&I_{2}&D_{6}^{*}&&&{xz(x+z)(x+t)}\\
 \cline{2-6}
 &1&\infty&0&&&{\rho(E_{w})=1}\\
\hline 
S_{7}&I_{2}&I_{2}&I_{2}&I_{2}&I_{4}&{x(x+t)(x+z-\lambda
  t)(x+z)}\\
 \cline{2-6}
 &0&1&\lambda&\lambda+1&\infty&{\rho(E_{w})=2}\\
\hline 
\end{array}
\]
For any of the above constructions it is not difficult to write down
also a Weierstrass equation.

\begin{remark}
The fibrations $S_{3}$ and $S_{4}$ have the same types of special
fibers but the quadruples of lines they are given by are not projectively equivalent.
However, they can be transported to each other by a Cremona transformation. 
\end{remark}

\subsection{Isogeny between $S_1$ and $S_2$}
\label{sec:isogeny}

Consider the elliptic fibration $S_{1}$ with congruence group
$\Gamma_{1}(4)\cap \Gamma(2)$. Computing the Euler and Hodge numbers of the
fiber product of $S_{1}$ and $S_{2}$ we obtain that the surfaces are isogenous.
For a fixed $t$ the fibers come from each other by doubling the lattice
in one direction. 

If we have any elliptic curve $E$ with equation of the type 
\[
y^2=x^3+Ax^2+Bx,
\]
then fixing the point at infinity as zero for the group structure the point
$e=(0,0)$ becomes a half period. Dividing $E$ by the map $E\ni p\mapsto p+e\in E$
we obtain as the quotient the curve 
\[
y^{2}=(x+A)(x^2-4B),
\]
and the quotient is given explicitly by the map 
\[
(x,y)\mapsto (x+\frac Bx, y-\frac B{x^2}).
\]

We will apply this to the elliptic fibrations $S_{1}$ and $S_{2}$. 
The surface $S_{1}$ can be given in local coordinates by the
(birationally equivalent) Weierstrass equation
\begin{equation}
  \label{eq:el2}
  y^{2}=x\left(x-(t^{2}-1)\right)\left(x-t^{2}\right).
\end{equation}

Let $S_{2}$ be the twist of $S_{1}$ by the automorphism of $\pz^1$
given by $t\mapsto\frac{t-1}{t+1}$. Similar equations for $S_{2}$ have the
form 
\begin{equation}
  \label{eq:el4}
  y^{2}=x\left(x-(t-1)^{2}\right)\left(x-(t+1)^{2}\right)
\end{equation}

Starting with \eqref{eq:el2} and replacing $x$ by $\frac{x}{4}$ and
$y$ by $\frac{y}{8}$ we obtain
\[
y^2 = x(x-4(t^2-1))(x-4t^2).
\]
Replacing $x$ by $x+2(t^2-1)$ we get
\[
y^2 = (x+2(t^2-1))(x-2(t^2-1))(x-2(t^2+1)).
\]
Now replacing $x$ by $x+\frac{(t^2-1)^2}{x}$ and $y$ by
$y(1-\frac{(t^2-1)^2}{x^2})$ we obtain \eqref{eq:el4}.

Composing the maps we find that the map $\phi:S_{2}\lra S_{1}$ given by
\begin{align*}
x \quad \mapsto & \quad \frac{x^2+16x(t^2-1)+80(t^2-1)^2}{4(x+8(t^2-1))}\\
y \quad \mapsto & \quad \frac{y}{8}\left( \frac{x^2 + 16x(t^2-1) + 48(t^2-1)}{(x + 8(t^2-1))^2}\right) 
\end{align*}
is a generically $2:1$ rational isogeny. 

\subsection{Isogeny between $S_{1}$ and $X_{1128}$}

Consider again the Weierstrass equation \eqref{eq:el2}
\begin{align*}
y^{2} & = x\left(x-(t^{2}-1)\right)\left(x-t^{2}\right)\\
& = x^3 + (1-2t^2)x^2 + t^2(t^2-1)x
\end{align*}
for the fibration $S_1$. Proceeding as in \ref{sec:isogeny}, i.e., setting
\[
(x,y) \mapsto \left(x+\frac{t^2(t^2-1)}{x}, y-\frac{t^2(t^2-1)}{x^2}\right),
\]
we obtain the extremal fibration given by
\[
y^2 = (x+1-2t^2)(x^2-4t^2(t^2-1))
\]
with singular fibers of type $I_1$, $I_1$, $I_2$, $I_8$. In \cite{MP} this
fibration is denoted by $X_{1128}$. The corresponding congruence group is
$\Gamma_0(8)\cap\Gamma_1(4)$.

Thus there is a generically $2:1$ rational map $\gamma : X_{1128}\lra S_1$.

\subsection{Pullback from $S_{5}$ to $S_{1}$}
\label{sec:pullback}

The fibrations $S_{1}\cong S_{2}$ and $S_{5}\cong S_{6}$ are extremal, in
\cite{MP} they are denoted by $X_{4422}$ and $X_{222}$. By \cite{MP} the
fibration $S_{1}$ can be obtained from $S_{5}$ by a base change. To see this
we substitute $t=t^{2}$ in the equation of $S_{5}$ in the above table
and we obtain an equation of the form 
\[
u^{2}=t^{2}x(x+1)(x+t^{2}).
\]
Substituting $x=x-t^{2}$ and taking a normalization we get equation~\eqref{eq:el2}.

Denote by $\psi$ the map $\psi:S_{1}\lra S_{5}$. It is a generically
$2:1$ rational map. There is also a similar map $\psi':S_{1}\lra S_{6}$.

\section{Rigid double octic Calabi--Yau threefolds}
\label{sec:rig}

In this section we list the fiber products of elliptic fibrations
which come from rigid double octic Calabi--Yau threefolds
constructed from arrangements of eight planes and listed in \cite{Meyer}
(we will also use the notations introduced there). Note that not all the
elliptic fibrations are semistable, but they all admit crepant
resolutions of singularities. We only have to consider the case of a
product of $I_{n}$ and $D_{m}^{*}$. Since $D_{m}^{*}$ contains double
lines the fiber product is singular along the product of such a line
and a node. After blowing up all double lines there will remain only
nodes which admit projective small resolutions.

The following table lists the double octics $X$ and their classification numbers
(as in the tables in \cite{Meyer}) and invariants and a description of the
corresponding (resolved) fiber products $Y$, including the types of the singular
fibers and the Euler and Hodge numbers. We also give the levels of the
corresponding weight four newforms. If the level is equal to 8 then we
add the types $T$, as introduced in section \ref{sec:list}.
\[
\def\arraystretch{1.3}
\begin{array}[t]{|c|c|c|c||c|c|c|c|c||c|}\hline
\multicolumn{4}{|c||}{\text{double octic } X}&
\multicolumn{5}{|c||}{\text{(resolved) fiber product } Y}&\\
\cline{1-9}
\rm No. &e&h^{11}&T&\text{singular fibers}&e&h^{11}&h^{12}&T&N\\\hline\hline
1&140&70&T_{70}^1&
\begin{array}[c]{ccc}
I_{2}&D_{6}^{*}&I_{2}\\
D_{6}^{*}&I_{2}&I_{2}
\end{array}&72&36&0&T_{36}&8\\\hline
3&124&62&&
\begin{array}[c]{cccc}
I_{4}&I_{4}&I_{2}&I_{2}\\
D_{6}^{*}&I_{2}&I_{2}&I_{0}
\end{array}&88&45&1&&32\\\hline
19&108&54&&
\begin{array}[c]{cccc}
I_{2}&I_{2}&I_{4}&I_{4}\\
I_{0}&D_{6}^{*}&I_{2}&I_{2}
\end{array}&64&33&1&&32\\\hline
32&100&50&T_{50}&
\begin{array}[c]{cccc}
I_{2}&I_{2}&I_{4}&I_{4}\\
I_{2}&I_{2}&D_{4}^{*}&I_{2}
\end{array}&80&40&0&T_{40}^2&8\\\hline
69&100&50&T_{50}&
\begin{array}[c]{cccc}
I_{2}&I_{2}&I_{4}&I_{4}\\
I_{2}&I_{2}&D_{4}^{*}&I_{2}
\end{array}&80&40&0&T_{40}^2&8\\\hline
93&92&46&T_{46}&
\begin{array}[c]{cccc}
I_{4}&I_{4}&I_{2}&I_{2}\\
I_{2}&I_{2}&D_{4}^{*}&I_{2}
\end{array}&64&32&0&T_{32}^2&8\\\hline
238&88&44&T_{44}&
\begin{array}[c]{cccc}
I_{2}&I_{2}&I_{4}&I_{4}\\
I_{2}&I_{2}&I_{4}&I_{4}
\end{array}&80&40&0&T_{40}^1&8\\\hline
239&80&40&&
\begin{array}[c]{ccccc}
I_{2}&I_{2}&I_{4}&I_{4}&I_{0}\\
I_{0}&I_{4}&I_{2}&I_{4}&I_{2}
\end{array}&64&34&2&&12\\\hline
240&80&40&&
\begin{array}[c]{ccccc}
I_{2}&I_{2}&I_{4}&I_{4}&I_{0}\\
I_{2}&I_{2}&I_{2}&I_{4}&I_{2}
\end{array}&64&33&1&&6\\\hline
241&80&40&T_{40}&
\begin{array}[c]{cccc}
I_{2}&I_{2}&I_{4}&I_{4}\\
I_{4}&I_{4}&I_{2}&I_{2}
\end{array}&64&32&0&T_{32}^1&8\\\hline
245&76&38&&
\begin{array}[c]{ccccc}
I_{2}&I_{2}&I_{4}&I_{4}&I_{0}\\
I_{4}&I_{2}&I_{2}&I_{2}&I_{2}
\end{array}&64&33&1&&6\\\hline

\end{array}
\]

Observe that the fiber products are rigid exactly for the cases where
the level is equal to 8. In the other cases, using the methods of \cite{Schuett, HulekVerrill}
we can identify codimension two modular motives in $H^{3}(Y)$. We can then show
the modularity of the remaining two-dimensional motive, either by counting
points over finite fields or by identifying it by the correspondence with
the modular rigid double octic.

Note also that the two level 32 newforms associated with the second and
the third example are different so these two Calabi--Yau threefolds can
not be in correspondence.

\section{Correspondences for level 8 rigid Calabi--Yau threefolds}
\label{sec:lev8}

In this section we use the results of sections \ref{sec:ref} and \ref{sec:rig}
to describe further correspondences between level 8 rigid Calabi--Yau threefolds.
In the end we compile all known correspondences into a picture.

Taking the fiber product of the maps $\phi$ and ${\rm id}_{S_{1}}$
we obtain a degree two correspondence between the Calabi--Yau threefolds that
are resolutions of the fiber products $S_{2}\times _{\pz^1}S_{1}$ and
$S_{1}\times _{\pz^1}S_{1}$ and hence a degree two correspondence between the
Calabi--Yau threefolds defined by arrangements no. 241 and 238, i.e.,
there are correspondences
\[
T_{32}^1 \stackrel{2:1}\longrightarrow T_{40}^1,\qquad
T_{40} \stackrel{4:2}\longrightarrow T_{44}.
\]

Taking the fiber product of the maps $\psi$ and $\psi'$ we obtain a
generically $4:1$ map from $S_{1}\times _{\pz^1}S_{1}$ to $S_{5}\times_{\pz^1}S_{6}$
and hence a degree four correspondence between the Calabi--Yau threefolds
defined by arrangements no. 238 and 1, i.e., there are correspondences
\[
T_{40}^1 \stackrel{4:1}\longrightarrow T_{36},\qquad
T_{44} \stackrel{8:2}\longrightarrow T_{70}^1.
\]

Taking the fiber product of the maps $\phi$ and ${\rm id}_{S_3}$ we obtain a
generically $4:1$ map from $S_{2}\times _{\pz^1}S_{3}$ to $S_{1}\times_{\pz^1}S_{3}$
and hence a degree two correspondence between the Calabi--Yau threefolds
defined by arrangements No. 32 and 93, i.e., there are correspondences
\[
T_{40}^2 \stackrel{2:1}\longrightarrow T_{32}^2,\qquad
T_{50} \stackrel{4:2}\longrightarrow T_{46}.
\]

Taking the self-fiber product of the map $\gamma$ we obtain a generically
$4:1$ map from $X_{1128}\times _{\pz^1}X_{1128}$ to $S_1\times _{\pz^1}S_1$
and hence a correspondence
\[
T_{70} \stackrel{4:1}\longrightarrow T_{40}^1.
\]

Finally consider the complete intersection $X$ from \ref{subsec:Nygaard} (type $T_{32}$)
and the Fermi threefold $Z$ from \ref{subsec:Fermi} (type $T_{70}$). J. Stienstra
constructed the $8:1$ rational map $X\lra Z$ given by
\[
x = \frac{y_0 + x_0}{y_0 - x_0}, \qquad
y = \frac{y_1 + x_1}{y_1 - x_1}, \qquad
z = \frac{y_2 + x_2}{y_2 - x_2}, \qquad
t = \frac{y_3 + x_3}{y_3 - x_3},
\]
i.e., there is a correspondence
\[
T_{32} \stackrel{8:1}\longrightarrow T_{70}.
\]
The map can be found in \cite[page 60]{NygaardvanGeemen} but there are some misprints.

\begin{remark}
We do not know whether the Calabi--Yau threefolds $T_{40}$, $T_{40}^1$,
$T_{40}^2$, $T_{40}^3$ (resp. $T_{32}$, $T_{32}^1$, $T_{32}^{2}$ and
$T_{70}$, $T_{70}^1$) are birational, however we belive that it is the case.
In general we expect that any two rigid Calabi--Yau threefolds with equal
Hodge numbers and the same $L$--series are birationally equivalent.
\end{remark}

The following picture contains known correspondences between level 8
rigid Calabi--Yau threefolds. To keep it concise we only included those
correspondences induced by explicit maps (including degrees).
Correspondences that were known before are marked with thicker lines.

\begin{center}
\includegraphics[width=9cm]{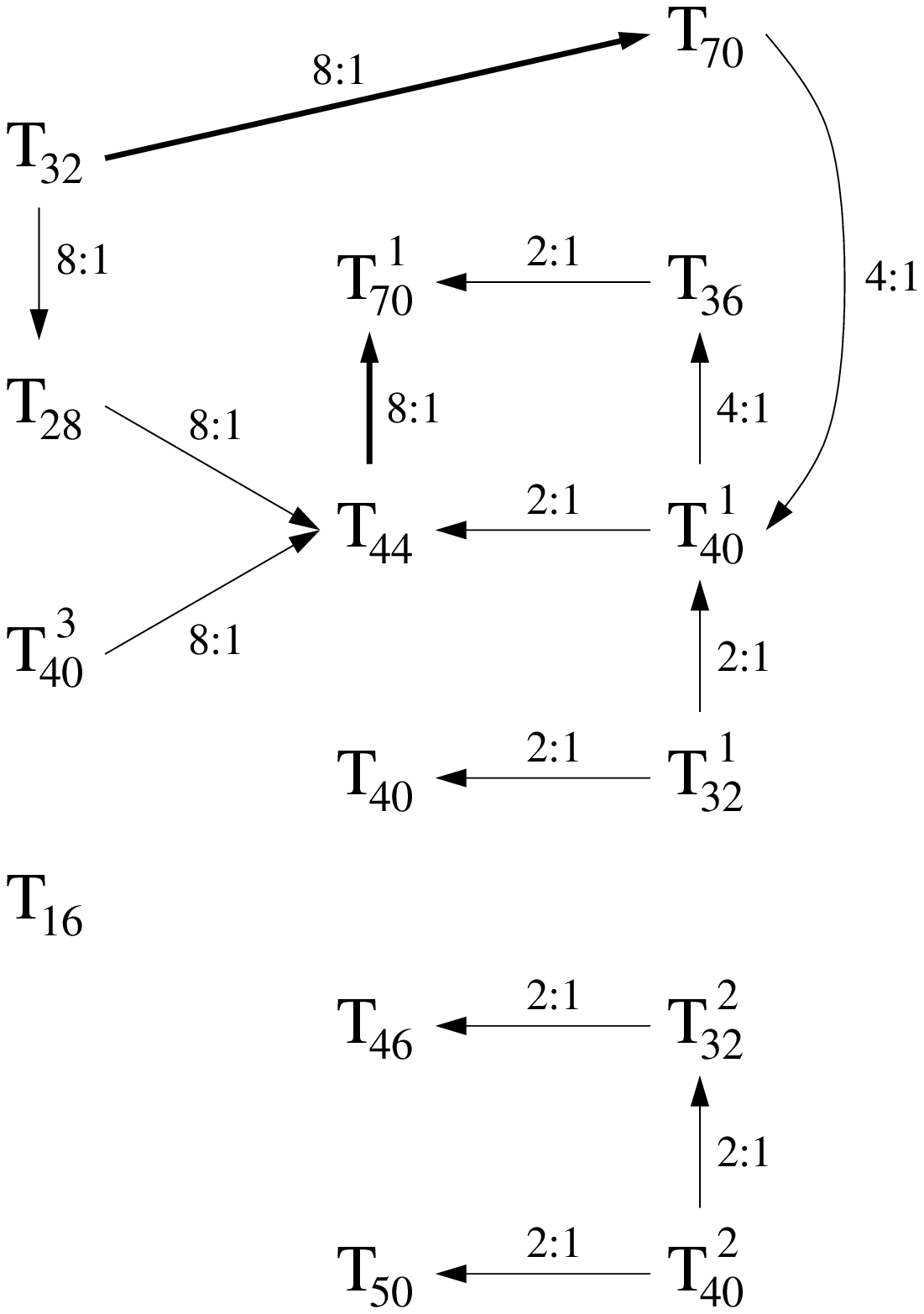}
\end{center}

Note that there is a correspondence between $T_{16}$ and a double covering
of $\pz^3$ branched along the union of two Kummer surfaces with 12 common
nodes. This correspondence is investigated in \cite{Meyer}. The Hodge
numbers of the double octic have not yet been computed but numerical
observations suggest that it is rigid.

\section{Level 6 double octic Calabi--Yau threefolds}
\label{sec:lev6}

In the table in section \ref{sec:rig} there are two level 6 rigid double octic
Calabi--Yau threefolds (constructed from arrangements no. 245 and 240).

Taking the fiber product of the maps $\phi$ and ${\rm id}_{S_7}$
we obtain a degree two correspondence between the Calabi--Yau threefolds that
are resolutions of the fiber products $S_{2}\times _{\pz^1}S_{7}$ and
$S_{1}\times _{\pz^1}S_{7}$ and hence a degree two correspondence between the
two double octics.

\subsection*{Acknowledgements}
The work on this paper was done during the first named author's stays at the
Institutes of Mathematics of the Johannes Gutenberg-Universit\"at Mainz
and the Universit\"at Hannover. He would like to thank both institutions
for their hospitality. The authors also would like to thank Prof. Duco van Straten,
Prof. Klaus Hulek and Matthias Sch\"utt for their help.

\end{document}